\documentclass[12pt]{article}
\usepackage{latexsym, amsbsy, amsfonts, amsmath, amsthm, textcomp, amssymb, fontenc}
\usepackage{graphics}
\usepackage[T1]{fontenc}
\usepackage{indentfirst}
\date{}

\usepackage{amscd}
\usepackage{amsfonts}
\usepackage{latexsym}
\usepackage{graphicx}
\usepackage{amsthm}
\usepackage{amsmath}

\def\rank{\mathrm{rank}}

\def \F{\mathcal{F}}

\newtheorem{thm}{Theorem}[section]

\newtheorem{lem}{Lemma}[section]

\theoremstyle{definition}

\newtheorem{rmk}{Remark} [section]

\title{An effective approach to Picard-Vessiot theory and the Jacobian Conjecture}

\begin{document}

\maketitle

\begin{center}
PAWE\L{} BOGDAN\\
Faculty of Mathematics and Computer Science, Jagiellonian University\\
ul. \L{}ojasiewicza 6, 30-348 Krak\'ow, Poland\\
e-mail: pawel.bogdan@uj.edu.pl
\end{center}

\begin{center}
ZBIGNIEW HAJTO\footnote{Z. Hajto acknowledges support of grant MTM2012-33830, Spanish Science Ministry.}\\
Faculty of Mathematics and Computer Science, Jagiellonian University\\
ul. \L{}ojasiewicza 6, 30-348 Krak\'ow, Poland\\
e-mail: zbigniew.hajto@uj.edu.pl
\end{center}

\begin{center}
 EL\.ZBIETA ADAMUS\footnote{E. Adamus acknowledges support of the Polish Ministry of Science and Higher Education.}\\
 Faculty of Applied Mathematics, AGH University of Science and Technology\\
al. Mickiewicza 30, 30-059 Krak\'ow, Poland\\
e-mail: esowa@agh.edu.pl
\end{center}

\vspace{20pt}

\begin{abstract}
In this paper we present a theorem concerning an equivalent statement of the
Jacobian Conjecture in terms of Picard-Vessiot extensions. Our theorem
completes the earlier work of T. Crespo and Z. Hajto which suggested an
effective criterion for detecting polynomial automorphisms of affine spaces.
We show a simplified criterion and give a bound on the number of wronskians
determinants which we need to consider in order to check if a given
polynomial mapping with non-zero constant Jacobian determinant is a
polynomial automorphism. Our method is specially efficient with cubic
homogeneous mappings introduced and studied in fundamental papers by H.
Bass, E. Connell, D. Wright and L. Dru\.{z}kowski.
\end{abstract}

\section{Introduction}

Let $K$ denote an algebraically closed field of characteristic zero. Let $%
n>0 $ be a fixed integer and let $F=(F_1, \ldots, F_n): K^n \rightarrow K^n$
be a polynomial mapping, i.e. $F_i \in K[X_1, \ldots, X_n]$ for $i=1,
\ldots, n$. We consider the Jacobian matrix $J_F = [\frac{\partial F_i}{%
\partial X_j}]_{1 \leq i,j \leq n}$. The Jacobian Conjecture states that if $%
\mathrm{det}(J_F)$ is a non-zero constant, then $F$ has an inverse, which is
also polynomial.

The Jacobian Conjecture is one of Stephen Smale's problems (cf. \cite%
{Sma}, Problem 16), which are a list of important problems in mathematics
for the twenty-first century. Originally the conjecture was formulated for $%
n=2$ by O. Keller (cf. \cite{Ke}). In 1982 H. Bass, E. Connell and D. Wright
(\cite{BCW}) showed that the general case follows from the case where $n\geq
2$ and $F=(X_{1}+H_{1},\ldots ,X_{n}+H_{n})$ and where each $H_{i}$ is zero
or homogeneous of degree 3. One year later L. Dru\.{z}kowski (\cite{D})
improved this result proving that if the Jacobian Conjecture is true for $%
n\geq 2$ and
\begin{equation}
F={\large (X_{1}+(\sum_{j=1}^{n}a_{1j}X_{j})^{3},\ldots
,X_{n}+(\sum_{j=1}^{n}a_{nj}X_{j})^{3}),}  \label{dru}
\end{equation}%
then it holds in general. A polynomial mapping $F$ of the form (\ref{dru})
with constant Jacobian is called a \emph{Dru\.{z}kowski mapping}. In 2001 Dru%
\.{z}kowski \cite{D2} proved that in his reduction (\ref{dru}) it is enough
to assume that the matrix $A=[a_{ij}]$ is nilpotent of degree 2, i.e. $%
A^{2}=0$.

In 2011 T. Crespo and Z. Hajto generalized a classical theorem of A.
Campbell (\cite{Ca}) by proving an equivalent statement of the Jacobian
Conjecture in terms of Picard-Vessiot extensions (cf. \cite{CH}, Theorem 2).
Condition 4 in Theorem 2 in the work of Crespo and Hajto suggested an
effective criterion for polynomial automorphisms of affine spaces. However,
the effectivity is obstructed by the big number of generalized wronskians
which have to be considered when the dimension of the affine space is
growing. In this paper we present a simplified criterion for a polynomial
automorphism of an affine space and prove that if the dimension of the space
is $n$ then it is enough to consider $\dfrac 1 2 n^2(n+1)-n$ generalized
wronskians. We believe that a deeper analysis of our algorithm may lead to
the proof of the Jacobian Conjecture.

Let $(\mathcal{\mathcal{F}}, \Delta_{\mathcal{F}})$ be a partial
differential field with an algebraically closed field of constants $C_{%
\mathcal{F}}$ and $\Delta_\F=\{\partial_{1},\ldots, \partial_{m}\}$. Let us
consider a \emph{linear partial differential system in matrix form} over $%
\mathcal{F}$, i.e. a system of equations of the form
\begin{equation}
\partial_{i}(Y)=A_iY, \, i=1, \ldots, m, \, A_i \in M_{n \times n}(\mathcal{F%
}).  \label{eq}
\end{equation}
A matrix $y \in GL_{n}(\mathcal{K})$, where $\mathcal{K}$ is a differential
field extension of $\mathcal{F}$, is called a \emph{fundamental matrix} for
the system (\ref{eq}) if $\partial_{i}(y)=A_iy$ for $i=1, \ldots, m$. We say
that the system (\ref{eq}) is \emph{integrable} if it has a fundamental
matrix. A differential field extension $(\mathcal{G}, \Delta _{\mathcal{G}})
$ of $(\mathcal{F}, \Delta_{\mathcal{F}})$ is a \emph{Picard-Vessiot
extension for the integrable system} (\ref{eq}) if the following holds: $C_{%
\mathcal{G}}=C_{\mathcal{F}}$, there exists a fundamental matrix $%
y=\{y_{ij}\} \in GL_n(\mathcal{G})$ and $\mathcal{G}$ is generated over $%
\mathcal{F}$ as a field by the entries of $y$, i.e. $\mathcal{G}=\mathcal{F}%
(\{y_{ij}\}_{1 \leq i,j, \leq n})$.

There is another definition of a Picard-Vessiot extension, formulated by
Kolchin in \cite{K}. Let $\mathcal{F}$ be a partial differential field of
characteristic zero with $\Delta_\F=\{\partial_{1},\ldots, \partial_{m}\}$
and algebraically closed field of constants $C_{\mathcal{F}}$. Let $\mathcal{%
G}$ be a differential field extension of $\mathcal{F}$. Let $Y_{1}, \ldots ,
Y_{n}$ denote indeterminates and let $\Theta$ denote the free commutative
multiplicative semigroup generated by the elements of $\Delta_{\F}$. So $\theta
\in \Theta$ is a differential operator of the form $\partial_1^{i_1}\ldots
\partial_{m}^{i_m}$, where $i_1, \ldots , i_m \in \mathbb{Z}_{+}\cup \{0\}$.
Let us denote by $\Theta(k)$ the subset of $\Theta$ of the elements of order
less than or equal to $k$. The determinant $\mathrm{det}(\theta_iy_j)_{1\leq
i,j \leq n}$ is called a \emph{generalized wronskian determinant} and denoted by $%
W_{\theta_1, \ldots \theta_n}(y_1, \ldots, y_n)$. Kolchin called $\mathcal{G}
$ a \emph{Picard-Vessiot extension of $\mathcal{F}$} if $C_{\mathcal{G}} =
C_{\mathcal{F}}$ and there exist $\eta_1, \ldots, \eta_n \in \mathcal{G} $
linearly independent over $C_{ \mathcal{F} }$ such that $\mathcal{G} =
\mathcal{F} \langle \eta_1, \ldots, \eta_n \rangle$ and
\begin{equation}
\forall \theta_1, \ldots \theta_n \in \Theta(n) : \, \frac{W_{\theta_1,
\ldots \theta_n}(\eta_1, \ldots, \eta_n)}{W_{\theta_{01}, \ldots
\theta_{0n}}(\eta_1, \ldots, \eta_n)} \in \mathcal{F}  \label{wr}
\end{equation}
for some fixed $\theta_{01}, \ldots \theta_{0n}$ such that $W_{\theta_{01},
\ldots \theta_{0n}}(\eta_1, \ldots, \eta_n) \neq 0$.

Theorem 1 in \cite{CH} establishes the equivalence between the two definitions
of Picard-Vessiot extension of partial differential fields presented above.
Theorem 2 in \cite{CH}, which is a differential version of the classical
theorem of Campbell, gives an equivalent formulation of the Jacobian
Conjecture. Let $K$ be an algebraically closed field of characteristic zero
and let $F=(F_1, \ldots, F_n):K^n \rightarrow K^n$ be a polynomial map such that $\mathrm{det}(J_F) = c \in K\setminus \{0\}$. We can equipp $%
K(x_1, \ldots, x_n)$ with the Nambu derivations, i.e. derivations $\delta_1,
\ldots, \delta_n$ given by
\begin{equation*}
\left(
\begin{array}{c}
\delta_1 \\
\vdots \\
\delta_n%
\end{array}
\right) = (J_F^{-1})^T \left(
\begin{array}{c}
\frac{\partial}{\partial x_1} \\
\vdots \\
\frac{\partial}{\partial x_1}%
\end{array}
\right).
\end{equation*}
Observe that $K\langle F_1, \ldots, F_n \rangle = K(F_1, \ldots, F_n)$, i.e.
$K(F_1, \ldots, F_n)$ is stable under $\delta_1, \ldots, \delta_n$. Moreover
if $\mathrm{det}(J_F) = 1$, then $J_F^{-1} = [\delta_jx_i]_{1 \leq i,j \leq
n}$.

The following theorem is a reformulation of theorems 1 and 2 in \cite{CH} in
the form we will use in the sequel.

\begin{thm}
Let $K$ and $F$ be as above. Then the following conditions are equivalent:

\begin{enumerate}
\item[1)] $F$ is a polynomial automorphism

\item[2)] The matrix

\begin{equation*}
W= \left[
\begin{array}{cccc}
1 & x_1 & \ldots & x_n \\
0 & \delta_1x_1 & \ldots & \delta_1x_n \\
\vdots & \vdots & \ddots & \vdots \\
0 & \delta_nx_1 & \ldots & \delta_nx_n%
\end{array}%
\right]
\end{equation*}

\noindent is a fundamental matrix for an integrable system

\begin{equation*}
\delta_k Y=A_k Y, \quad k=0, \ldots, n,
\end{equation*}

\noindent where we are taking $\delta_0=\mathrm{id}$, with $A_k \in
M_{(n+1)\times (n+1)}(K(F_1,\dots,F_n))$.
\end{enumerate}

\label{thm2}
\end{thm}

\section{Wronskian criterion}

Theorem \ref{thm2} gives a method of checking if a given polynomial map $F$
is a polynomial automorphism. If we denote $x_0=1$, then we may write $%
W=[\delta_ix_j]_{i,j=0,1,\ldots,n}$. Let us assume that $det W=1$ (which is
equivalent to $\mathrm{det}(J_F) = 1$). We are going to find $A_k=\delta_kW
\cdot W^{-1} $, $k=1,2,\ldots, n$ in order to check if the entries of $A_k$%
's lie in $K(F_1,\dots,F_n)$. We have that

\begin{equation*}
\delta_kW=\left[%
\begin{array}{cccc}
0 & \delta_kx_1 & \ldots & \delta_kx_n \\
0 & \delta_k\delta_1x_1 & \ldots & \delta_k\delta_1x_n \\
0 & \delta_k\delta_2x_1 & \ldots & \delta_k\delta_2x_n \\
\ldots &  &  & \ldots \\
0 & \delta_k\delta_ix_1 & \ldots & \delta_k\delta_ix_n \\
\ldots &  &  & \ldots \\
0 & \delta_k\delta_nx_1 & \ldots & \delta_k\delta_nx_n%
\end{array}
\right] = \big[ \omega^k_{ij} \big] _{i,j=0,1, \ldots, n}, \quad \text{where}
\quad \omega^k_{ij}=\delta_k\delta_ix_j.
\end{equation*}
Let us find $W^{-1}=\big([D_{ij}]_{i,j=0,1, \ldots, n} \big)^T$, where $%
D_{ij}$ denote the adjoint determinant of the element $\delta_ix_j$ of
matrix $W$. We obtain that
\begin{equation*}
D_{00}=1\quad \mathrm{and} \quad \forall j\geq 1: \, D_{0j}=0,
\end{equation*}

\begin{equation*}
D_{i0}=(-1)^{i+1+1} \left|
\begin{array}{ccc}
x_1 & \ldots & x_n \\
\delta_1x_1 & \ldots & \delta_1 x_n \\
\ldots & \ldots & \ldots \\
\delta_{i-1}x_1 & \ldots & \delta_{i-1} x_n \\
\delta_{i+1}x_1 & \ldots & \delta_{i+1} x_n \\
\ldots & \ldots & \ldots \\
\delta_nx_1 & \ldots & \delta_nx_n%
\end{array}
\right| = (-1)^{i+0} \mathrm{det} \big( [\delta_sx_t]_{s=0,1,\ldots, n; \, s
\neq i; \, t=1, \ldots n} \big).
\end{equation*}
For $i,j \geq 1$ we get

\begin{equation*}
D_{ij}=(-1)^{i+j} \left|
\begin{array}{ccccccc}
1 & x_1 & \ldots & x_{j-1} & x_{j+1} & \ldots & x_n \\
0 & \delta_1x_1 & \ldots & \delta_1x_{j-1} & \delta_1x_{j+1} & \ldots &
\delta_1x_n \\
\vdots & \vdots &  & \vdots & \vdots &  & \vdots \\
0 & \delta_{i-1}x_1 & \ldots & \delta_{i-1}x_{j-1} & \delta_{i-1}x_{j+1} &
\ldots & \delta_{i-1}x_n \\
0 & \delta_{i+1}x_1 & \ldots & \delta_{i+1}x_{j-1} & \delta_{i+1}x_{j+1} &
\ldots & \delta_{i+1}x_n \\
\vdots & \vdots &  & \vdots & \vdots &  & \vdots \\
0 & \delta_nx_1 & \ldots & \delta_nx_{j-1} & \delta_nx_{j+1} & \ldots &
\delta_nx_n%
\end{array}
\right| .
\end{equation*}
So $D_{ij}=(-1)^{i+j} \mathrm{det} \big( [\delta_sx_t]_{s,t=\mathbf{%
0,1,\ldots,n}; \, s\neq i, t \neq j} \big)
=(-1)^{i+j} \mathrm{det} \big( [\delta_sx_t]_{s,t=\mathbf{1,\ldots,n}; \,
s\neq i, t \neq j} \big)$ and consequently $W^{-1}=[B_{ij}]_{i,j=0,1,
\ldots, n}$, where
\begin{equation*}
B_{ij}=(-1)^{i+j}D_{ji} = (-1)^{i+j} \left|
\begin{array}{ccccccc}
1 & x_1 & \ldots & x_{i-1} & x_{i+1} & \ldots & x_n \\
0 & \delta_1x_1 & \ldots & \delta_1x_{i-1} & \delta_1x_{i+1} & \ldots &
\delta_1x_n \\
\vdots & \vdots &  & \vdots & \vdots &  & \vdots \\
0 & \delta_{j-1}x_1 & \ldots & \delta_{j-1}x_{i-1} & \delta_{j-1}x_{i+1} &
\ldots & \delta_{j-1}x_n \\
0 & \delta_{j+1}x_1 & \ldots & \delta_{j+1}x_{i-1} & \delta_{j+1}x_{i+1} &
\ldots & \delta_{j+1}x_n \\
\vdots & \vdots &  & \vdots & \vdots &  & \vdots \\
0 & \delta_nx_1 & \ldots & \delta_nx_{i-1} & \delta_nx_{i+1} & \ldots &
\delta_nx_n%
\end{array}
\right| .
\end{equation*}
We compute $A_k=\big[ a^k_{ij} \big] _{i,j=0,1, \ldots, n} = \delta_kW \cdot
W^{-1}$. We obtain $a^k_{i0}=0$, i.e. the first column (i.e. the column
indexed by j=0) is a zero column. Moreover for $j \geq 1$

\begin{equation*}
a^k_{0j} = \sum_{r=1}^n \delta_k\delta_0x_r \cdot B_{rj}=\sum_{r=1}^n
\delta_kx_r \cdot B_{rj}= \delta_kx_1 \cdot (-1)^{1+j} \left|
\begin{array}{ccc}
\delta_1x_2 & \ldots & \delta_1x_n \\
\ldots & \ldots & \ldots \\
\delta_{j-1}x_2 & \ldots & \delta_{j-1} x_n \\
\delta_{j+1}x_2 & \ldots & \delta_{j+1} x_n \\
\ldots & \ldots & \ldots \\
\delta_nx_2 & \ldots & \delta_nx_n%
\end{array}
\right| + \ldots
\end{equation*}
\begin{equation*}
\ldots+ \delta_kx_n \cdot (-1)^{n+j} \left|
\begin{array}{ccc}
\delta_1x_1 & \ldots & \delta_1x_{n-1} \\
\ldots & \ldots & \ldots \\
\delta_{j-1}x_1 & \ldots & \delta_{j-1} x_{n-1} \\
\delta_{j+1}x_1 & \ldots & \delta_{j+1} x_{n-1} \\
\ldots & \ldots & \ldots \\
\delta_nx_1 & \ldots & \delta_nx_{n-1}%
\end{array}
\right| = \left|
\begin{array}{cccc}
\delta_1x_1 & \delta_1x_2 & \ldots & \delta_1x_n \\
\vdots & \vdots & \ddots & \vdots \\
\delta_{j-1}x_1 & \delta_{j-1}x_2 & \ldots & \delta_{j-1} x_n \\
\delta_kx_1 & \delta_kx_2 & \ldots & \delta_kx_n \\
\delta_{j+1}x_1 & \delta_{j+1} x_2 & \ldots & \delta_{j+1} x_n \\
\vdots & \vdots & \ddots & \vdots \\
\delta_nx_1 & \delta_nx_2 & \ldots & \delta_nx_n%
\end{array}
\right| = \left\{
\begin{array}{ccc}
0 & ; & k \neq j \\
1 & ; & k=j%
\end{array}
\right. .
\end{equation*}
If $i,j \geq 1$, then $a^k_{ij}=\sum_{r=1}^n \delta_k\delta_ix_r \cdot
B_{rj} $, this means we have
\begin{equation*}
a^k_{ij} = \delta_k\delta_ix_1 \cdot (-1)^{1+j} \left|
\begin{array}{ccc}
\delta_1x_2 & \ldots & \delta_1x_n \\
\ldots & \ldots & \ldots \\
\delta_{j-1}x_2 & \ldots & \delta_{j-1} x_n \\
\delta_{j+1}x_2 & \ldots & \delta_{j+1} x_n \\
\ldots & \ldots & \ldots \\
\delta_nx_2 & \ldots & \delta_nx_n%
\end{array}
\right| + \ldots+ \delta_k\delta_ix_n \cdot (-1)^{n+j} \left|
\begin{array}{ccc}
\delta_1x_1 & \ldots & \delta_1x_{n-1} \\
\ldots & \ldots & \ldots \\
\delta_{j-1}x_1 & \ldots & \delta_{j-1} x_{n-1} \\
\delta_{j+1}x_1 & \ldots & \delta_{j+1} x_{n-1} \\
\ldots & \ldots & \ldots \\
\delta_nx_1 & \ldots & \delta_nx_{n-1}%
\end{array}
\right| =
\end{equation*}
\begin{equation*}
=\left|
\begin{array}{cccc}
\delta_1x_1 & \delta_1x_2 & \ldots & \delta_1x_n \\
\vdots & \vdots & \ddots & \vdots \\
\delta_{j-1}x_1 & \delta_{j-1}x_2 & \ldots & \delta_{j-1} x_n \\
\delta_k\delta_ix_1 & \delta_k\delta_ix_2 & \ldots & \delta_k\delta_ix_n \\
\delta_{j+1}x_1 & \delta_{j+1} x_2 & \ldots & \delta_{j+1} x_n \\
\vdots & \vdots & \ddots & \vdots \\
\delta_nx_1 & \delta_nx_2 & \ldots & \delta_nx_n%
\end{array}
\right|
\end{equation*}

The total number of considered determinants is $n(n+1)^2$, since for every $%
\delta_k$ we have $(n+1)^2$ of them and $k=1, \ldots, n$. However for each $%
A_k$ we can ignore the first row and the first column (i.e. the row and the
column indexed by 0), since they consist of constant elements. Consequently,
we can omit $2n+1$ of elements for each $A_k$. So there are $n^3$ wronskians
left. We can easily observe that for every $j=1, \ldots, n$ and for $k\neq i$
we have $a^k_{ij}=a^i_{kj}$. So we can omit $\binom{n}{2}$ of determinants
for each $j$. Due to the lemma given below we can omit even more
determinants.

\begin{lem}
Let $(K,^{\prime})$ be a differential field and let $A
=[a_{ij}]_{i,j=1,\ldots, n} \in GL_n(K)$ be a nonsingular matrix with
entries in $K$. Then
\begin{equation*}
(\mathrm{det}{A})^{\prime}= \left|
\begin{array}{cccc}
a_{11} & a_{12} & \ldots & a_{1n} \\
a_{21} & a_{22} & \ldots & a_{2n} \\
\vdots & \vdots & \ddots & \vdots \\
a_{n1} & a_{n2} & \ldots & a_{nn}%
\end{array}
\right| ^{\prime}=
\end{equation*}
\begin{equation*}
=\left|
\begin{array}{cccc}
a_{11}^{\prime} & a_{12}^{\prime} & \ldots & a_{1n}^{\prime} \\
a_{21} & a_{22} & \ldots & a_{2n} \\
\vdots & \vdots & \ddots & \vdots \\
a_{n1} & a_{n2} & \ldots & a_{nn}%
\end{array}
\right|+ \left|
\begin{array}{cccc}
a_{11} & a_{12} & \ldots & a_{1n} \\
a_{21}^{\prime} & a_{22}^{\prime} & \ldots & a_{2n}^{\prime} \\
\vdots & \vdots & \ddots & \vdots \\
a_{n1} & a_{n2} & \ldots & a_{nn}%
\end{array}
\right|+\ldots + \left|
\begin{array}{cccc}
a_{11} & a_{12} & \ldots & a_{1n} \\
a_{21} & a_{22} & \ldots & a_{2n} \\
\vdots & \vdots & \ddots & \vdots \\
a_{n1}^{\prime} & a_{n2}^{\prime} & \ldots & a_{nn}^{\prime}%
\end{array}
\right|=
\end{equation*}
\begin{equation*}
=\left|
\begin{array}{cccc}
a_{11}^{\prime} & a_{12} & \ldots & a_{1n} \\
a_{21}^{\prime} & a_{22} & \ldots & a_{2n} \\
\vdots & \vdots & \ddots & \vdots \\
a_{n1}^{\prime} & a_{n2} & \ldots & a_{nn}%
\end{array}
\right|+ \left|
\begin{array}{cccc}
a_{11} & a_{12}^{\prime} & \ldots & a_{1n} \\
a_{21} & a_{22}^{\prime} & \ldots & a_{2n} \\
\vdots & \vdots & \ddots & \vdots \\
a_{n1} & a_{n2}^{\prime} & \ldots & a_{nn}%
\end{array}
\right|+\ldots + \left|
\begin{array}{cccc}
a_{11} & a_{12} & \ldots & a_{1n}^{\prime} \\
a_{21} & a_{22} & \ldots & a_{2n}^{\prime} \\
\vdots & \vdots & \ddots & \vdots \\
a_{n1} & a_{n2} & \ldots & a_{nn}^{\prime}%
\end{array}
\right|
\end{equation*}
\label{detlem}
\end{lem}

Let us use lemma (\ref{detlem}) to differentiate $\mathrm{det} W=1$ with
respect to each $\delta_k$, $k=1, \ldots, n$. We get that
\begin{equation*}
\delta_1 \left|
\begin{array}{ccc}
\delta_1x_1 & \ldots & \delta_1x_n \\
\vdots & \ddots & \vdots \\
\delta_nx_1 & \ldots & \delta_nx_n%
\end{array}%
\right| = \left|
\begin{array}{ccc}
\delta_1^2x_1 & \ldots & \delta_1^2x_n \\
\vdots & \ddots & \vdots \\
\delta_nx_1 & \ldots & \delta_nx_n%
\end{array}%
\right| + \left|
\begin{array}{ccc}
\delta_1x_1 & \ldots & \delta_1x_n \\
\delta_1\delta_2x_1 & \ldots & \delta_1\delta_2x_n \\
\delta_3x_1 & \ldots & \delta_3x_n \\
\vdots & \ddots & \vdots \\
\delta_nx_1 & \ldots & \delta_nx_n%
\end{array}%
\right| +\ldots+ \left|
\begin{array}{ccc}
\delta_1x_1 & \ldots & \delta_1x_n \\
\vdots & \ddots & \vdots \\
\delta_{n-1}x_1 & \ldots & \delta_{n-1}x_n \\
\delta_1\delta_nx_1 & \ldots & \delta_1\delta_nx_n%
\end{array}%
\right| =0  \label{delta1}
\end{equation*}

\begin{equation*}
\ldots \qquad \ldots \qquad \ldots
\end{equation*}

\begin{equation*}
\delta_n \left|
\begin{array}{ccc}
\delta_1x_1 & \ldots & \delta_1x_n \\
\vdots & \ddots & \vdots \\
\delta_nx_1 & \ldots & \delta_nx_n%
\end{array}%
\right| = \left|
\begin{array}{ccc}
\delta_n\delta_1x_1 & \ldots & \delta_n\delta_1x_n \\
\delta_2x_1 & \ldots & \delta_2x_n \\
\vdots & \ddots & \vdots \\
\delta_nx_1 & \ldots & \delta_nx_n%
\end{array}%
\right| + \left|
\begin{array}{ccc}
\delta_1x_1 & \ldots & \delta_1x_n \\
\delta_n\delta_2x_1 & \ldots & \delta_n\delta_2x_n \\
\delta_3x_1 & \ldots & \delta_3x_n \\
\vdots & \ddots & \vdots \\
\delta_nx_1 & \ldots & \delta_nx_n%
\end{array}%
\right| +\ldots+ \left|
\begin{array}{ccc}
\delta_1x_1 & \ldots & \delta_1x_n \\
\vdots & \ddots & \vdots \\
\delta_{n-1}x_1 & \ldots & \delta_{n-1}x_n \\
\delta_n^2x_1 & \ldots & \delta_n^2x_n%
\end{array}%
\right| =0  \label{deltan}
\end{equation*}

Let us go back to considerations concerning the matrix $A_k$. We can use the
equations given above to observe that for each $k=1, \ldots, n$ we have $%
a^k_{11} + \ldots + a^k_{nn}=0$. So for example
\begin{equation*}
a^k_{kk}=-a^k_{11}-\ldots-a^k_{k-1,k-1}-a^k_{k+1,k+1}- \ldots - a^k_{nn} .
\end{equation*}
So we can omit $n$ determinants more. Hence it is enough to check the
following number of wronskian determinants
\begin{equation}  \label{eq:n}
n^3 - n \cdot \binom{n}{2}-n=n^3-\frac{1}{2}n^2(n-1)-n=\frac{1}{2}n^2(n+1)-n.
\end{equation}

\noindent Let us observe that the number given in (\ref{eq:n}) is optimal,
e.g. for $n=2$, we have to consider 4 wronskian determinants.

\section{Examples}

In this section in order to explain how our criterion works for detecting
polynomial automorphisms we shall present two explicit examples.\newline

\emph{Example 1.} Let us consider a well-known wild automorphism: the Nagata
automorphism:\newline

\noindent $%
F_{1}=x_{1}-2x_{2}(x_{3}x_{1}+x_{2}^{2})-x_{3}(x_{3}x_{1}+x_{2}^{2})^{2} $%
\newline
$F_{2}=x_{2}+x_{3}(x_{3}x_{1}+x_{2}^{2})$\newline
$F_{3}=x_{3}$\newline

\noindent Using the computer algebra system Maple18 we first compute that $%
\mathrm{det}(J_{F})=1$ and next the entries of the matrices $\big[a_{ij}^{k}%
\big]_{i,j=1,2,3}$, for $k=1,2,3$. We obtain the following results:\newline

\textbf{k=1:}\newline
\noindent $%
a_{11}^{1}=-a_{22}^1-a_{33}^1=-2x_{3}^{3}(-2x_{1}x_{3}^{3}-2x_{2}^{2}x_{3}^{2}-2x_{2}x_{3}+1)=4F_2F_3^4-2F_3^3
$\newline

\noindent $a_{12}^{1}=-2x_{3}^{5}=-2F_{3}^{5}$\newline

\noindent$a_{13}^{1}=0$\newline

\noindent $%
a_{21}^{1}=(-4x_{1}x_{3}^{4}-4x_{2}^{2}x_{3}^{3}-4x_{2}x_{3}^{2}+2x_{3})(-2x_{1}x_{3}^{3}-2x_{2}^{2}x_{3}^{2} -2x_{2}x_{3}+1)=8F_2^2F_3^3-8F_2F_3^2+2F_3
$\newline

\noindent$%
a_{22}^{1}=(-4x_{1}x_{3}^{4}-4x_{2}^{2}x_{3}^{3}-4x_{2}x_{3}^{2}+2x_{3})x_{3}^{2}= -4F_2F_3^4+2F_3^3
$\newline

\noindent $a_{23}^{1}=0$\newline

\noindent$%
a_{31}^{1}=-4x_1^3x_3^8-12x_1^2x_2^2x_3^7-12x_1x_2^4x_3^6-4x_2^6x_3^5-12x_1^2x_2x_3^6-24x_1x_2^3x_3^5 -12x_2^5x_3^4+10x_1^2x_3^5+8x_1x_2^2x_3^4-2x_2^4x_3^3+16x_1x_2x_3^3+12x_2^3x_3^2+6x_2^2x_3+2x_2= 4F_2^2F_3+4F_1F_2F_3^3-2F_1F_3^2+2F_2
$\newline

\noindent $%
a_{32}^{1}=2x_1^2x_3^7+4x_1x_2^2x_3^6+2x_2^4x_3^5+4x_1x_2x_3^5+4x_2^3x_3^4-4x_1x_3^4-2x_2^2x_3^3-2x_2x_3^2-2x_3= -2F_1F_3^4-2F_2F_3^2-2F_3
$\newline

\noindent$a_{33}^{1}=0$\newline

\textbf{k=2:}\newline
\noindent $a_{11}^{2}=a_{21}^{1}$\newline
\noindent $a_{12}^{2}=a_{22}^{1}$\newline
\noindent $a_{13}^{2}=a_{23}^{1}=0$\newline

\noindent $a_{21}^{2}=
16x_1^3x_3^8+48x_1^2x_2^2x_3^7+48x_1x_2^4x_3^6+16x_2^6x_3^5+48x_1^2x_2x_3^6+96x_1x_2^3x_3^5 +48x_2^5x_3^4-24x_1^2x_3^5+24x_2^4x_3^3-48x_1x_2x_3^3-32x_2^3x_3^2+12x_1x_3^2-12x_2^2x_3+12x_2 = 16F_2^3F_3^2-24F_2^2F_3+12F_2
$\newline

\noindent $%
a_{22}^{2}=-a_{11}^2-a_{33}^2=-8x_1^2x_3^7-16x_1x_2^2x_3^6-8x_2^4x_3^5-16x_1x_2x_3^5-16x_2^3x_3^4+8x_1x_3^4+8x_2x_3^2-2x_3= -8F_2^2F_3^3+8F_2F_3^2-2F_3
$\newline

\noindent $a_{23}^{2}=0$\newline

\noindent $a_{31}^{2}=
-8x_1^4x_3^9-32x_1^3x_2^2x_3^8-48x_1^2x_2^4x_3^7-32x_1x_2^6x_3^6-8x_2^8x_3^5
-32x_1^3x_2x_3^7-96x_1^2x_2^3x_3^6-96x_1x_2^5x_3^5-32x_2^7x_3^4+24x_1^3x_3^6+24x_1^2x_2^2x_3^5 -24x_1x_2^4x_3^4-24x_2^6x_3^3+64x_1^2x_2x_3^4+96x_1x_2^3x_3^3+32x_2^5x_3^2-10x_1^2x_3^3+36x_1x_2^2x_3^2 +38x_2^4x_3-12x_1x_2x_3+4x_2^3+2x_1= 8F_1F_2^2F_3^2-8F_1F_2F_3+8F_2^3+2F_1
$\newline

\noindent $a_{32}^{2}=
4x_1^3x_3^8+12x_1^2x_2^2x_3^7+12x_1x_2^4x_3^6+4x_2^6x_3^5+12x_1^2x_2x_3^6
+24x_1x_2^3x_3^5+12x_2^5x_3^4-10x_1^2x_3^5-8x_1x_2^2x_3^4+2x_2^4x_3^3-16x_1x_2x_3^3-12x_2^3x_3^2-6x_2^2x_3-2x_2= -4F_1F_2F_3^3+2F_1F_3^2-4F_2^2F_3-2F_2
$\newline

\noindent $a_{33}^{2}=0$\newline

\textbf{k=3:}\newline
\noindent $a_{11}^{3}=a_{31}^1$\newline
\noindent $a_{12}^{3}=a_{32}^{1}$\newline
\noindent $a_{13}^{3}=a_{33}^{1}=0$\newline
\noindent $a_{21}^{3}=a_{31}^2$\newline
\noindent $a_{22}^{3}=a_{32}^2$\newline
\noindent $a_{23}^{3}=a_{33}^{2}=0$\newline

\noindent $a_{31}^{3}=
4x_1^5x_3^{10}+20x_1^4x_2^2x_3^9+40x_1^3x_2^4x_3^8+40x_1^2x_2^6x_3^7+20x_1x_2^8x_3^6 +4x_2^{10}x_3^5+20x_1^4x_2x_3^8+80x_1^3x_2^3x_3^7+120x_1^2x_2^5x_3^6+80x_1x_2^7x_3^5+20x_2^9x_3^4-18x_1^4x_3^7 -32x_1^3x_2^2x_3^6+12x_1^2x_2^4x_3^5+48x_1x_2^6x_3^4+22x_2^8x_3^3-64x_1^3x_2x_3^5-152x_1^2x_2^3x_3^4 -112x_1x_2^5x_3^3-24x_2^7x_3^2+16x_1^3x_3^4-36x_1^2x_2^2x_3^3-100x_1x_2^4x_3^2-48x_2^6x_3+28x_1^2x_2x_3^2 +8x_1x_2^3x_3-16x_2^5-2x_1^2x_3+8x_1x_2^2= 4F_1^2F_2F_3^2-2F_1^2F_3+8F_1F_2^2
$\newline

\noindent $a_{32}^{3}=
-2x_1^4x_3^9-8x_1^3x_2^2x_3^8-12x_1^2x_2^4x_3^7-8x_1x_2^6x_3^6-2x_2^8x_3^5-8x_1^3x_2x_3^7-24x_1^2x_2^3x_3^6-24x_1x_2^5x_3^5 -8x_2^7x_3^4+8x_1^3x_3^6+12x_1^2x_2^2x_3^5-4x_2^6x_3^3+20x_1^2x_2x_3^4+32x_1x_2^3x_3^3+12x_2^5x_3^2-4x_1^2x_3^3+8x_1x_2^2x_3^2 +10x_2^4x_3+4x_2^3-2x_1 = -2F_1^2F_3^3-4F_1F_2F_3-2F_1
$\newline

\noindent $a_{33}^{3}=-a_{11}^3-a_{22}^3=0$\newline

\emph{Example 2.} Recently Dan Yan (\cite{Y}) has proved that the Jacobian
Conjecture is true for the Dru\.{z}kowski mappings in dimension $n\leq 9$,
however only in the case when the matrix $A$ (cf.(\ref{dru})) has no zeros on its
diagonal, and for general $n$ and $\rank A \leq 4$.  Moreover Michiel de Bondt in his thesis \cite{dB} proved the validity of the Jacobian Conjecture for all Dru\.{z}kowski mappings in dimension $n\leq 8$. Let us consider the following Dru\.{z}kowski mapping in dimension 13.

\bigskip

$$\begin{array}{lll}
F_{ 1 } & = & X_{ 1 } + \left( \frac{1}{6} X_{4} + \frac{1}{6} X_{5} - \frac{1}{3} X_{6} - \frac{1}{6} X_{7} - \frac{1}{6} X_{8} + \frac{1}{3} X_{9} + X_{13} \right)^3 \\
F_{ 2 } & = & X_{ 2 } + \left( \frac{1}{6} X_{4} + \frac{1}{6} X_{5} - \frac{1}{3} X_{6} - \frac{1}{6} X_{7} - \frac{1}{6} X_{8} + \frac{1}{3} X_{9} -  X_{13} \right)^3 \\
F_{ 3 } & = & X_{ 3 } + \left( \frac{1}{6} X_{4} + \frac{1}{6} X_{5} - \frac{1}{3} X_{6} - \frac{1}{6} X_{7} - \frac{1}{6} X_{8} + \frac{1}{3} X_{9} \right)^3 \\
F_{ 4 } & = & X_{ 4 } + \left( \frac{1}{6} X_{1} + \frac{1}{6} X_{2} - \frac{1}{3} X_{3} + X_{12} \right)^3 \\
F_{ 5 } & = & X_{ 5 } + \left( \frac{1}{6} X_{1} + \frac{1}{6} X_{2} - \frac{1}{3} X_{3} -  X_{12} \right)^3 \\
F_{ 6 } & = & X_{ 6 } + \left( \frac{1}{6} X_{1} + \frac{1}{6} X_{2} - \frac{1}{3} X_{3} \right)^3 \\
F_{ 7 } & = & X_{ 7 } + \left( -\frac{1}{3} X_{3} + \frac{1}{6} X_{10} + \frac{1}{6} X_{11} + X_{13} \right)^3 \\
F_{ 8 } & = & X_{ 8 } + \left( -\frac{1}{3} X_{3} + \frac{1}{6} X_{10} + \frac{1}{6} X_{11} -  X_{13} \right)^3 \\
F_{ 9 } & = & X_{ 9 } + \left( -\frac{1}{3} X_{3} + \frac{1}{6} X_{10} + \frac{1}{6} X_{11} \right)^3 \\
F_{ 10 } & = & X_{ 10 } + \left( \frac{1}{6} X_{4} + \frac{1}{6} X_{5} - \frac{1}{3} X_{6} - \frac{1}{6} X_{7} - \frac{1}{6} X_{8} + \frac{1}{3} X_{9} + X_{12} \right)^3 \\
F_{ 11 } & = & X_{ 11 } + \left( \frac{1}{6} X_{4} + \frac{1}{6} X_{5} - \frac{1}{3} X_{6} - \frac{1}{6} X_{7} - \frac{1}{6} X_{8} + \frac{1}{3} X_{9} -  X_{12} \right)^3 \\
F_{ 12 } & = & X_{ 12 }   \\
F_{ 13 } & = & X_{ 13 }
\end{array}
$$

In the above example $rank(A)=5$.
The computation of wronskians is involved,
therefore we have presented it separately in our website $http://crypto.ii.uj.edu.pl/galois/$.

\begin{rmk}
In his landmark paper \cite{Ca} L.A. Campbell studied in fact general covering
maps. Let us observe that our computational approach can be used as well for
detecting Galois coverings. In this case we can not assume that the Jacobian
determinant is a non-zero constant, however the computations are
analogous.
\end{rmk}


\begin{thebibliography}{99}
\bibitem{BCW} H. Bass, E. Connell, D. Wright, \emph{The Jacobian
Conjecture: Reduction of Degree and Formal Expansion of the Inverse},
Bulletin of the American Mathematical Society 7 (1982), 287-330.

\bibitem{dB} M. de Bondt, \emph{Homogeneous Keller maps}, Ph. D. thesis, July 2007,
http://webdoc.ubn.ru.nl/mono/b/bondt$_{-}$m$_{-}$de/homokema.pdf

\bibitem{Ca} L.A. Campbell, \emph{A condition for a polynomial map to be
invertible}, Math. Annalen 205 (1973), 243-248.

\bibitem{CH} T. Crespo, Z. Hajto, \emph{Picard-Vessiot theory and the
Jacobian problem}, Israel Journal of Mathematics \textbf{186}, 2011, pp.
401-406.

\bibitem{D} L. M. Dru\.zkowski, \emph{An Effective Approach to Keller's
Jacobian Conjecture}, Math. Ann. 264 (1983), 303-313.

\bibitem{D2} L. M. Dru\.zkowski, \emph{New reduction in the Jacobian
conjecture}, Univ. Iagell. Acta Math. 39 (2001), 203-206.

\bibitem{Ke} O.H. Keller, \emph{Ganze Cremona Transformationen},
Monatsh. Math. Phys. 47 (1939), 299-306.

\bibitem{K} E. R. Kolchin, \emph{Picard-Vessiot theory of partial
differential fields}, Proceedings of the American Mathematical Society
\textbf{3}, 1952, pp. 596-603.

\bibitem{Sma} S. Smale, \emph{Mathematical Problems for the Next Century}%
, Mathematical Intelligencer, 20 (1998), 7-15.

\bibitem{Y} D. Yan, \emph{A note on the Jacobian Conjecture}, Linear
Algebra and its Applications 435, 2110-2113 (2011).
\end{thebibliography}
\end{document}